\documentclass[11pt]{article}
\usepackage{amssymb,amsmath,latexsym}
\usepackage[mathscr]{eucal}

\newtheorem{thm}{Theorem}[section]
\newtheorem{cor}[thm]{Corollary}
\newtheorem{lem}[thm]{Lemma}
\newtheorem{prop}[thm]{Proposition}

\newtheorem{remark}[thm]{Remark}
\numberwithin{equation}{section}

\def\FF{{\cal F}}
\def\cS{{\cal S}}
\def\K{{\cal K}}
\def\E{{\mathbb E}}
\def\P{{\mathbb P}}
\def\HH{{\mathbb H}}
\def\RR{{\mathbb R}}
\def\CC{{\mathbb C}}
\def\DD{{\mathbb D}}

\def\A{{\mathbb A}}

\def\dd{{\cal D}}

\def\pf{\noindent{\bf Proof.} }
\def\qed{{\hfill $\Box$ \bigskip}}
\def\wh{\widehat}
\def\wt{\widetilde}
\def\dis{\displaystyle}
\def\eps{\varepsilon}

\def\c{\check}

\def\1{{\bf 1}}

\def\x{{\bf x}}
\def\y{{\bf y}}
\def\s{{\bf s}}
\def\SS{{\bf S}}
\def\H{{\bf H}}

\def\<{\langle}
\def\>{\rangle}

\begin{document}
\title{\bf Stochastic Komatu-Loewner evolutions and SLEs}
\author{{\bf Zhen-Qing Chen},
\ {\bf Masatoshi Fukushima}
 \ and \ {\bf Hiroyuki Suzuki}}  
 \date{(April 27, 2016)}
\maketitle

\begin{abstract}
Let $D={\mathbb H}\setminus \bigcup_{j=1}^N C_j$ be a standard slit domain, where $\HH$ is the upper half plane and $C_j,1\le j\le N,$ are mutually disjoint horizontal line segments in $\HH$.  A stochastic Komatu-Loewner evolution denoted by ${\rm SKLE}_{\alpha,b}$ has been introduced in \cite{CF} as a family $\{F_t\}$ of random growing hulls with $F_t\subset D$ driven by a diffusion  process $\xi(t)$ on $\partial {\mathbb H}$ that is determined by certain continuous homogeneous functions $\alpha$ and $b$ defined on the space $\cS$ of all labelled standard slit domains.  We aim at identifying the distribution of a suitably reparametrized ${\rm SKLE}_{\alpha,b}$ with that of the Loewner evolution on ${\mathbb H}$ driven by the path of a certain continuous semimartingale and thereby relating the former to the distribution of ${\rm SLE}_{\alpha^2}$ when $\alpha$ is a constant.  We then prove that, when $\alpha$ is a constant, ${\rm SKLE}_{\alpha,b}$ up to some random hitting time and modulo a time change has the same distribution as ${\rm SLE}_{\alpha^2}$ under a suitable Girsanov transformation. We further show that 
a reparametrized ${\rm SKLE}_{\sqrt{6},-b_{\rm BMD}}$ has the same distribution  as ${\rm SLE}_6$, where $b_{\rm BMD}$ is the BMD-domain constant indicating the discrepancy of $D$ from $\HH$ relative to Brownian motion with darning (BMD in abbreviation).  A key ingredient of the proof is a hitting time analysis for the absorbing Brownian motion on ${\mathbb H}.$   We also revisit and examine the locality property of 
${\rm SLE}_6$ in several canonical domains.  Finally K-L equations and SKLEs for other canonical multiply connected planar domains than the standard slit one are recalled and examined.
 \end{abstract}

\medskip

\smallskip\noindent 
{\bf AMS 2010 Mathematics Subject Classification}: Primary 60J67, 60J70; Secondary 30C20, 60H10, 60H30

\smallskip\noindent
{\bf Keywords and phrases:} standard slit domain, stochastic Komatu-Loewner evolution, SLE, absorbing Brownian motion, locality

\section{Introduction}
A subset $A$ of the upper half-plane $\HH$
is called an $\HH$-{\it hull} if $A$ is bounded closed in $\HH$ and $\HH\setminus A$ is simply connected.  Given an $\HH$-hull $A$, there exists a unique conformal map $f$ (one-to-one analytic function) from $\HH\setminus A$ onto $\HH$ satisfying a hydrodynamic normalization
at infinity
\[ f(z)=z+\frac{a}{z}+o(1/|z|),\ z\to\infty.\]
Such a map will be called a {\it canonical Riemann map from} $\HH\setminus A$
and the constant $a$ is called the {\it half-plane capacity} of $A$ relative to $f$.

We consider a simple ODE called a chordal {\it Loewner differential equation}
\begin{equation}\label{1.1}
\frac{d z(t)}{dt}= -2\pi \Psi^\HH(z(t),\xi(t)),\quad z(0)=z\in \HH, 
\end{equation} 
where
\[ \Psi^\HH(z,\xi)=-\frac{1}{\pi} \frac{1}{z-\xi},\ z\in \HH,\ \xi \in \partial \HH,\]
that is the so-called {\it complex Poisson kernel} for the absorbing Brownina motion on $\HH$ because
\[
\Im \Psi^\HH(z,\xi)=\frac1\pi \frac {y}{(x-\xi)^2+y^2} \quad \hbox{for } z=x+iy \in \HH
\hbox{ and } \xi \in \partial \HH ,
\]
is the classical Poisson kernel in upper half space $\HH$.

Given a continuous function $\xi(t)\in \HH,\  0\le t<\infty,$ the Cauchy problem of the ODE \eqref{e:1} admits a unique solution $z(t),\ 0\le t<t_z,$ with the maximal interval of definition  $[0,t_z)$.   Define
\begin{equation}
K_t=\{z\in \HH: t_z\le t\},\quad t\ge 0.
\end{equation}
Then $K_t$ is an $\HH$-hull and $z(t)$ is the canonical Riemann map from $\HH\setminus K_t.$  The family $\{K_t: t\ge 0\}$ of growing hulls is called the chordal {\it Loewner evolution driven by} $\xi(t),\ 0\le t<\infty.$

Let $B(t)$ be   one-dimensional Brownian motion and $\kappa$ be a positive constant.  The random Loewner evolution driven by the sample path of $B(\kappa t),\ 0\le t<\infty,$ is called the {\it stochastic Loewner evolution} and is denoted by ${\rm SLE}_\kappa.$      It was introduced by Oded Schramm \cite{S} in his consideration of critical two-dimensional lattice models in statistical physics and their scaling limits.  
It is  now also called the {\it Schramm-Loewner evolution}.  Remarkable features as the {\it locality property} of ${\rm SLE}_6$ and the {\it restriction property} of ${\rm SLE}_{8/3}$ were then revealed (\cite{LSW1, LSW2}).  ${\rm SLE}_\kappa$ was shown in \cite{RS} to be generated by a continuous curve in the sense that, there exists a continuous path $\gamma: [0,\infty) \;\mapsto \;\overline \HH$ such that $\HH\setminus K_t$ is identical with the unbounded connected component of $\HH\setminus \gamma[0,t]$ a.s. for each $t>0,$ and furthermore $\gamma$ was shown to be a 
simple curve when $\kappa\le 4$, self-intersecting when $4<\kappa \leq 8$ and space-filling when $\kappa >8$.            

Several attempts have been made to extend both of the Loewner equation and the associated SLE from simply connected planar domains to multiply connected ones.
    Recently, motivated by \cite{BF1, BF2, L2}, Chen-Fukushima-Rohde \cite{CFR} and Chen-Fukushima \cite{CF}  studied 
   Komatu-Loewner equation and   stochastic Komatu-Loewner evolution, respectively,
in standard slit domains of finite multiplicity.  
 Stochastic Komatu-Loewner evolution, denoted by ${\rm SKLE}_{\alpha,b}$, 
 is a family of conformal maps   that are determined by two functions $\alpha$
 and $b$ on the slit space $S$ to be described below.  They generate  an increasing family of random growing $\HH$-hulls. 
 
 The main purpose of this paper is to study the geometry of 
${\rm SKLE}_{\alpha,b}$-hulls.
We show that,  after   a suitably reparametrization,  ${\rm SKLE}_{\alpha,b}$-hulls have the same disturbution
as  that of the Loewner evolution on $\HH$ driven by a continuous semi-martingale.
In particular, we show that when function $\alpha$ is a constant, 
 after a reparametrization and under an equivalent martingale measure, 
  ${\rm SLE}_{\alpha , b}$ has the same distribution as the chordal ${\rm SLE}_{\alpha^2}$ in $\HH$ up to a stopping time.
 Hence when $\alpha$ is a positive constant, we  conclude that  ${\rm SLE}_{\alpha , b}$-hulls are generated by continuous paths which
  are simple if $\alpha \leq 2$, self-intersecting  if $2<\alpha \leq 2\sqrt{2}$ and space-filling when $\alpha >2\sqrt{2}$. 
 
Fix $N\geq 1$. A {\it standard slit domain} (of $N$ slits) is a domain of the type $D=\HH \setminus \bigcup_{j=1}^N C_j$ for mutually disjoint line segments $C_j\subset \HH$ parallel to $\partial\HH.$   
The collection of all labelled standard slit domains (of $N$ slits) is denoted by $\dd.$
For $D\in \dd$,  let $z_k=x_k+iy_k, \ z_k^r=x_k^r+iy_k$ be the left and right endpoints of the $k$th slit $C_k$ of $D$.
It is characterized by $\y :=(y_1, \dots, y_N)$, $\x:=(x_1, \dots, x_N)$ and $\x^r:= (x_1^r, \dots, x_N^r)$
with the property that $\y>0$, $\x <\x^r$,  and either $x_j^r<x_k$ or $x_k^r <x_j$ whenever $y_j=y_k$ for $j\not= k$. 
Here for   vectors $\x, \y\in \RR^N$, $\y >0$ means each coordinate is strictly larger than 0; and 
$\x <\y$ means $\y - \x >0$.  With this characterization, the space $\dd$ can be identified with the following 
subset of $\RR^{3N}$  
\begin{eqnarray*}
 \cS=\left\{(\y,\x,\x')\in \RR^{3N}: \y> {\bf 0},\ \x<\x',     \ {\rm either}\ x_j'<x_k\ {\rm or}\ x_k'<x_j\ {\rm whenever}\ y_j=y_k,\ j\neq k 
  \right\}. 
\end{eqnarray*}
    For $\s\in \cS,$ denote by $D(\s)$ the corresponding element in $\dd$.  For $\xi\in \RR,$ we denote by $\wh \xi\in \RR^{3N}$ the $3N$-vector whose first $N$-components are equal to $0$ and the rest are equal to $\xi$.

For $\s\in \cS,$ we denote by $\Psi_\s(z,\xi),\; z\in D(\s),\; \xi\in \partial \HH,$ the {\it BMD-complex Poisson kernel} for $D=D(\s)$, namely, the unique analytic function in $z\in D$ vanishing at $\infty$ whose imaginary part is the Poisson kernel 
for the {\it Brownian motion with darning} (BMD) for $D$ (see \cite{CFR}).  

A function $f$ on $\cS$ is called {\it homogeneous} with degree $0$ (resp. $-1$) if $f(c\s)=f(\s)$ (resp. $f(c\s)=c^{-1} f(\s)$) for 
every positive constant  $c>0.$  It is said to satisfy the {\it local Lipschitz condition} if the following property holds:

\smallskip\noindent 
{\bf (L)}\ For any $\s \in \cS $ and any finite open interval $J\subset \RR,$ there exist a neighborhood $U(\s )$ of $\s $ in $\cS $ and a constant $L>0$ such that
\begin{equation}\label{e:4.1}
|f(\s^{(1)}-\wh \xi)-f(\s^{(2)}-\wh \xi)|\le L\;|\s^{(1)}-\s^{(2)}|
\quad \hbox{for } \s^{(1)},  \s^{(2)}\in U(\s ) \hbox{ and } \xi \in J.
\end{equation}

We consider the strong solution $(\xi(t), \s(t))\in \partial \HH\times \cS$ of the following stochastic differential equation (SDE) 
\begin{equation}\label{SDE1}
\begin{cases}
\xi(t)=\xi+\int_0^t\alpha(\s(t)-\wh\xi(t))dB_s+\int_0^t b(\s(t)-\wh\xi(t))ds \\
\s_j(t)=\s_j+\int_0^t  b_j(\s(t)-\wh\xi(t))ds,\quad t\ge 0,\quad 1\le j\le 3N, 
\end{cases}
\end{equation}
where $\alpha(\s)$ (resp. $b(\s)$) is a   
  homogneous function    on $\cS$ of degree $0$ (resp. $-1$) satisfying    
condition {\bf (L)} and  
\begin{equation}\label{SDE2}
 b_j(\s) :=
\begin{cases}
-2\pi\Im \Psi_\s(z_j, 0),\quad & 1\le j\le N,\\
-2\pi\Re \Psi_\s(z_j, 0),\quad & N+1\le j\le 2N,\\
-2\pi\Re \Psi_\s(z_j', 0),\quad & 2N+1\le j\le 3N.
\end{cases}
\end{equation}
It is known (see \cite{CF}) that $b_j(\s)$ is  a homogeneous function on $\cS$ of degree $-1$ satisfying condition {\bf (L)}.

Putting  the solution $(\xi(t), \s(t))$ of \eqref{SDE1} into the {\it Komatu-Loewner equation} introduced in \cite{CFR}, we consider 
the equation 
\begin{equation}\label{KLE}
\frac{d}{dt}g_t(z)=-2\pi \Psi_{\s(t)}(g_t(z), \xi(t)) \  \hbox{  with }   g_0(z)=z\in D . 
\end{equation}
The above equation has  a unique maximal solution $g_t(z),\ t\in [0,t_z),$ passing through  $G=\bigcup_{t\in [0,\zeta)}\{t\}\times  D_t,$ where $D_t=D(\s(t))$ and $D=D_0.$  Define 
\begin{equation}\label{hull}
F_t=\{z\in D:t_z\le t\},\quad t\ge 0.
\end{equation}

For $D\in \dd$ and an $\HH$-hull $A\subset D,$ the conformal map $f$ from $D\setminus A$ onto another set in $\dd$ satisfying the hydrodynamic normalization at infinity will be called the 
 {\it canonical map from} 
$D\setminus A$.
The set $F_t$ defined by \eqref{hull} is an $\HH$-hull and $g_t$ is the canonical map 
 from 
$D\setminus F_t.$  This family of growing hulls $\{F_t\}$ is denoted by ${\rm SKLE}_{\alpha,b}$ and will be called a {\it stochastic Komatu-Loewner evolution}.
${\rm SLE}_\kappa$ can be viewed as a special case of ${\rm SKLE}_{\alpha,b}$ where no slit is present, $\alpha$ is constant with $\alpha^2=\kappa$ and $b=0.$

For ${\rm SKLE}_{\alpha, b}$-hull $F_t$ defined by \eqref{hull}, we can  consider the canonical Riemann map $g_t^0(z)$ from $\HH\setminus F_t,$ the half-plane capacity $a(t)$ of $F_t$ relative to $g_t^0$ and a reparametrization $\{\check F_t\}$ of $\{F_t\}$ defined by $\c F_t=F_{a^{-1}(2t)},\; t\ge 0.$
 With  ${\rm SKLE}_{\alpha,b}$  reparametrized in this way, it is shown in  Theorem \ref{T:4.1} of this paper that it has the same distribution
 as the Schramm-Loewner evolution in $\HH$ driven by a continuous semimartingale $\c U(t).$  We then prove that, when $\alpha$ is a constant, ${\rm SKLE}_{\alpha,b}$ up to some random hitting time and modulo a time change,  has the same distribution as ${\rm SLE}_{\alpha^2}$,
  under a suitable Girsanov transformation; see Theorem \ref{T:4.3}. 
  Moreover, we show in Theorem \ref{T:4.2} that   ${\rm SKLE}_{\sqrt{6}, -b_{\rm BMD}}$, after a reparametrization,
    has the same distribution as ${\rm SLE}_6$, where $b_{\rm BMD}$ is the BMD-domain constant defined by \eqref{e:11b} that describes the discrepancy of a standard slit domain from $\HH$ relative to BMD. 

In order to establish Theorem \ref{T:4.1}  with rigor, we need to show that

\medskip\noindent
{\bf (C)}\ $g_t^0(z)$ is jointly continuous in $(t,z)\in [0,a]\times (\overline \HH\setminus F_a)$ for each $a>0.$

\medskip\noindent
A proof of this property will be carried out in Section 3 by combining the probabilistic representation
  of $\Im g_t^0(z)$ and $\Im g_t(z)$ obtained in \cite{CFR} in terms of the absorbing Brownian motion $Z^\HH$ on $\HH$ and BMD for $D$ 
with the continuity of $g_t(z)$ in $t$ 
that is the solution of the ODE \eqref{KLE}.  A key ingredient of the proof is a hitting time analysis for $Z^\HH.$
 
It is established in  \cite[Theorem 6.11]{CF}
that ${\rm SKLE}_{\sqrt{6}, -b_{\rm BMD}}$ 
enjoys a  locality property.  
 In relation to this and the present Theorem 4.2, 
we will present in Section 5  a first rigorous proof of the locality of the chordal ${\rm SLE}_6$ 
   in the sense of \cite{LSW3},  and 
 point out the missing pieces or gaps in other locality proofs in literature.

 In the final Section 6, we recall and examine Komatu-Loewner equations and stochastic Komatu-Loewner evolutions for other canonical multiply connected domains than the standard slit one.
 
\section{Riemann maps $\{g_t^0\}$ and a process $U(t)$ associated with SKLE}

Let $\alpha >0$ and $b$ be   homogeneous functions  on $\cS$ of degree $0$ and $-1$, respectively, 
that are local Lipschitz continuous.
Let    $(\xi(t), \s(t)),\; t<\zeta,$ be the strong solution of the associated SDE \eqref{SDE1}  and  
 $\{F_t\}$ be ${\rm SKLE}_{\alpha,b}$, namely, the family of growing hulls \eqref{hull} on  $D=D(\s(0))=\HH\setminus K,\ K=\cup_{j=1}^N C_j,$  
driven by $(\xi(t), \s(t)).$ 

Denote by  $g_t$ the canonical map from $D\setminus F_t$ onto $D_t=D(\s(t))$, $\Phi$  the identity map from $D$ into $\HH$,
 and $g_t^0$  the canonical Riemann map from $\HH\setminus F_t$ onto $\HH$.
According to \cite[Theorem 5.8]{CF}, $\{F_t\}$ is right continuous with limit $\xi(t)$ in the sense that 
\begin{equation}\label{e:1}
\bigcap_{\eps>0} \overline{g_t(F_{t+\eps}\setminus F_t)}=\xi(t).
\end{equation}  

Define
\begin{equation}\label{e:2}
\Phi_t(z)=g_t^0\circ\Phi\circ g_t^{-1}(z)  \qquad \hbox{for } z\in D_t=D(\s(t)).
\end{equation}

\begin{lem}
$\Phi_t$ admits an analytic extension to $D_t\cup \Pi D_t\cup\partial \HH$ by the Schwarz reflection.  Here $\Pi z=\overline z, \ z\in \HH.$ 
\end{lem}

\pf\ Take an arbitrary smooth Jordan arc $\Gamma$ in $\HH$ with two end points $z_1, z_2\in \partial \HH$ such that the open region $V$ enclosed by $\Gamma$ and the line segment connecting $z_1, z_2$ contains the set $F_t$ with $\overline V\cap K=\emptyset.$  
Clearly, $V_t:=g_t(V)$ is the 
open region enclosed by $g_t(\Gamma)$ and the line segment $\ell_t$ connecting $g_t(z_i)$, $i=1,2$.   In view of \eqref{e:1}, $\xi(t)$ is located 
in the interior of the line segment $\ell_t$. 
 Furthermore, $\Phi_t$ is a Riemann map from the Jordan domain $V_t$ onto 
the Jordan domain $g^0_t(V)$, which is 
enclosed by $g_t^0(\Gamma)$ and the line segment $\ell_t^0$ connecting $g_t^0(z_i)$,  $i=1,2$,  
and $\Phi_t$ maps 
 $\ell_t$ onto $\ell_t^0$ homeomorphically.  Thus $\Phi_t$ admits a
Schwarz reflection.  \qed

Define        
\begin{equation}\label{e:3}
 U(t)=\Phi_t(\xi(t)).
\end{equation}
We then have
\begin{equation}\label{e:4}
\bigcap_{\eps>0} \overline{g_t^0(F_{t+\eps}\setminus F_t)}=U(t),
\end{equation}  
because, by \eqref{e:1} and \eqref{e:2},
\[\bigcap_{\eps>0} \overline{g_t^0(F_{t+\eps}\setminus F_t)}
=\bigcap_{\eps>0} \overline{g_t^0\circ\Phi(F_{t+\eps}\setminus F_t)}
=\bigcap_{\eps>0} \overline{\Phi_t\circ g_t(F_{t+\eps}\setminus F_t)}=\Phi_t(\xi(t))=U(t).\]
For $D\in \dd$ and for an $\HH$-hull $A\subset D$, 
we denote by ${\rm Cap}^\HH(A)$ (resp. ${\rm Cap}^D(A)$) the half-plane capacity of $A$ relative to the canonical Riemann map $g_A^\HH$   
from $\HH\setminus A$ (resp. the canonical map $g_A^D$ from $D\setminus A$).  
\[{\rm Cap}^\HH(A)=\lim_{z\to\infty} z(g_t^\HH(z)-z),\quad {\rm Cap}^D(A)=\lim_{z\to\infty} z(g_t^D(z)-z).\]
Set 
$a(t) :={\rm Cap}^\HH(F_t)$ and $b(t) :={\rm Cap}^D(F_t).$

\begin{lem}\label{L:2.2} The right derivative of $a(t)$ 
\begin{equation}\label{e:5}
\frac{d^+ a(t)}{dt}:= \lim_{{\partial} \downarrow 0} \frac{ a(t+{\partial} )-a (t)}{{\partial}} =  
 2\Phi_t'(\xi(t))^2,
\end{equation}
\end{lem}

\pf\ For a set $A\subset \HH,$ we put ${\rm rad}(A)=\sup_{z\in A}|z|.$
For a fixed $t>0,$ let $K_\eps=g_t(F_{t+\eps}\setminus F_t).\; \;\eps>0.$  By \cite[Theorem 5.8 (iii)]{CF}, ${\rm rad}(K_\eps-\xi(t)) \to 0$ as $\eps \to 0.$  Hence by  the capacity comparison theorem 
\cite[Theorem 7.1]{CF}, we have
\begin{equation}\label{e:6}
{\rm Cap}^\HH(K_\eps)-{\rm Cap}^{D_t}(K_\eps)=o(\eps),\quad \eps\to 0.
\end{equation}
On the other hand, by  \cite[(3.8)]{L1}, 
\begin{equation}\label{e:7}
 a(t+\eps)-a(t)={\rm Cap}^\HH(g_t^0(F_{t+\eps}\setminus F_t)).
\end{equation}
Since 
$g_t^0(F_{t+\eps}\setminus F_t)=g_t^0\circ\Phi(F_{t+\eps}\setminus F_t)
=\Phi_t\circ g_t(F_{t+\eps}\setminus F_t)=\Phi_t(K_\eps),$
we obtain from \cite[(4.15)]{L1}, \eqref{e:6} and \eqref{e:7}
\[ a(t+\eps)-a(t)={\rm Cap}^\HH(\Phi_t(K_\eps))
=\Phi'_t(\xi(t))^2 {\rm Cap}^\HH(K_\eps)+o(\eps)
=\Phi'_t(\xi(t))^2 {\rm Cap}^{D_t}(K_\eps)+o(\eps),\]
which can be seen in an analogous manner to \eqref{e:7} to be equal to
\[\Phi'_t(\xi(t))^2 ({\rm Cap}^D(F_{t+\eps})-{\rm Cap}^D(F_t))+o(\eps)
=\Phi'_t(\xi(t))^2 (b(t+\eps)-b(t))+o(\eps).
\]
As $b(t)=2t$ by \cite[Theorem 5.12]{CF}, we arrive at \eqref{e:5}.
\qed

\begin{prop}\label{P:2.3}  It holds that
\begin{equation}\label{e:8}
\frac{d^+ g_t^0(z)}{dt}= \frac{2\Phi_t'(\xi(t))^2}{g_t^0(z)-U(t)},\quad z\in \HH\setminus F_t.
\end{equation}
in the right derivative sense.
\end{prop}
\pf  Denote by $Q$ the family of all $\HH$-hulls.  
According to \cite[p69, Propositon 3.46]{L1},
\begin{equation}\label{e:9}
g_{A-x}^\HH(z)=g_A^\HH(z+x)-x,\quad {\rm Cap^\HH}(A-x)={\rm Cap^\HH}(A) \quad A\in Q,\ x\in \RR,
\end{equation}
and, there exists a constant $c>0$ such that, for any $A\in Q$ and any $z$ with $|z|\ge 2{\rm rad}(A),$
\begin{equation}\label{e:10}
\left| z-g_A^\HH(z)+\frac{{\rm Cap}^\HH(A)}{z}\right| \le c\frac{{\rm rad}(A){\rm Cap}^\HH(A)}{|z|^2}.
\end{equation}

For $z\in \HH\setminus F_s,$ we get from \eqref{e:7}, \eqref{e:9} and \eqref{e:10}
\begin{eqnarray*}
 &&  g_{s+\eps}^0(z)-g_s^0(z) \\
&=& g^\HH_{g_s^0(F_{s+\eps}\setminus F_s)}(g_s^0(z))-g_s^0(z)\\
&  = & g^\HH_{g_s^0(F_{s+\eps}\setminus F_s)-U(s)}(g_s^0(z)-U(s))-(g_s^0(z)-U(s)) \\
&=&\frac{a(s+\eps)-a(s)}{g_s^0(z)-U(s)} 
+{\rm rad}(g_s^0(F_{s+\eps}\setminus F_s)-U(s))(a(s+\eps)-a(s))
O(1/(g_s^0(z)-U(s))^2)).
\end{eqnarray*}
The formula  \eqref{e:8} now follows from \eqref{e:4} and \eqref{e:5}.
\qed

  To show that the right derivative in  
Proposition \ref{P:2.3} can be strengthened to true derivative, 
  we need the following proposition, whose proof   is postponed to  next section.
 
\begin{prop}\label{P:2.4}
\ The Riemann maps 
$\{g_t^0\}$ enjoys the property {\rm \bf (C)} stated in Section 1.
\end{prop}

In the rest of this section, we shall take the validity of this proposition for granted.
The following lemma can then be shown exactly in the same way as the proof of \cite[Proposition 6.7 (i)]{CF}.  

\begin{lem}\label{L:2.5}
$\Phi_t(z),\ \Phi_t'(z),\ \Phi_t''(z)$  are jointly continuous in $(t,z)\in [0,\zeta)\times (D_t\cup \partial\HH).$
\end{lem}

By the property {\bf (C)} 
and the above lemma, the right hand side of \eqref{e:8} becomes continuous in $t$ and so \cite[Lemma 4.3]{L1} applies in getting the following thoerem.

\begin{thm}\label{T:2.6} \ 
 $g_t^0(z)$ is continuously differentiable in $t$ and {\rm\eqref{e:8}} becomes a genuine ODE:
\begin{equation}\label{e:11}
\frac{d g_t^0(z)}{dt}= \frac{2\Phi_t'(\xi(t))^2}{g_t^0(z)-U(t)},\quad  z\in \HH\setminus F_t.
\end{equation}
\end{thm}

\begin{remark}\label{R:2.7} \rm Strengthening from right time derivative in Proposition \ref{P:2.3} to 
the genuine time derivative in Theorem \ref{T:2.6} is very important 
since \eqref{e:8} does not uniquely characterize the conformal maps
$\{g^0_t(z)\}$. 
 This is because while the solution to 
\eqref{e:11} is unique, equation \eqref{e:8} may have numerous solutions. To see this, consider the case that $K=\emptyset$,
that is,  upper half space $\HH$ with no slits.  In this case, $\Phi_t (z)=z$ and \eqref{e:11} is the chordal Loewner equation
with driving function $U(t)$. So for each $z\in \HH$, 
\begin{equation}\label{e:2.12} 
\frac{d  g^0_t (z) }{dt}= \frac{2}{g^0_t (z) -U (t)}, \quad z (t) =z,
\end{equation}
 has a unique continuous solution $g^0_t(z)$ up to time
$t_z$ when $g^0_t$ and $U(t)$ collide. However, equation 
\begin{equation}\label{e:2.13} 
\frac{d^+ z (t) }{dt}= \frac{2}{z(t) -U (t)}, \quad z (t) =z,
\end{equation}
has infinitely many solutions. For instance, take any $\eps \in (0, \zeta_z)$ and
define $z(t)= g^0_t (z)$ for $t\in (0, \eps]$. Let $z(\eps)$ be any value in $\HH$.
Let  $\wt g^0_t (z(\eps))$, $0\leq t< t_{z(\eps)}$  be the unique solution of
$$
\frac{d  \wt g^0_t (z(\eps) ) }{dt}= \frac{2}{\wt g^0_t (z (\eps)) -U (t+\eps)}, \quad \wt g^0_0( z(\eps ))
 =z (\eps) . 
 $$
Define $z(t)= \wt g^0_{t-\eps} (z(\eps))$ for $t\in [\eps, \eps +  t_{z(\eps)} )$.   
Then $\{z(t); 0\leq t < \eps +  t_{z(\eps)}\}$ 
is a solution to equation  \eqref{e:2.13}.  
Indeed, we see by \cite[Lemma 4.3]{L1} that the solution $z(t)$ of \eqref{e:2.13} coincides with the solution $g_t^0(z)$ of \eqref{e:2.12} 
if and only if $z(t)$ is (left) continuous.  
\qed 
\end{remark}

For $\s \in \cS$, let $b_{\rm BMD}(\s)$ be the BMD-domain constant for the slit domain $D (\s)$ introduced in \cite[\S 6.1]{CF}:
\begin{equation}\label{e:11b}
b_{\rm BMD}(\s)=2\pi \lim_{z\to 0}\left(\Psi_\s(z,0)+\frac{1}{\pi z}\right).
\end{equation}

\begin{thm} \label{T:2.7}
\ The process $U(t)$ on $\partial\HH$ admits a semi-martingale decomposition
\begin{eqnarray}\label{e:12}
dU(t)&=& \Phi_t'(\xi(t))\alpha(\s(t)-\wh\xi(t))dB_t
 + \Phi_t'(\xi(t))\left( b_{\rm BMD}(\s(t)-\wh\xi(t))+b(\s(t)-\wh\xi(t))\right) dt\nonumber\\
&& + \Phi_t''(\xi(t))\left( -3+\frac12\alpha(\s(t)-\wh\xi(t))^2\right) dt.
\end{eqnarray}
\end{thm}  

\pf\ For a differentiable function $f_t(z):=f(t, z)$ defined on on open subset of $\RR_+ \times \CC$, 
we will use $\dot f$ and $f'$ to denote its partial derivative in $t$
and in $z\in \CC$,  respectively.  
Let $f_t (z) =g_t^{-1}$ (z).  Then 
\[\dot f_t(z)=2\pi f_t'(z)\Psi_{\s(t)}(z,\xi(t)),\quad z\in D_t,\]
 and $\Phi_t=g_t^0\circ\Phi\circ f_t$ by \eqref{e:2}.  
Thus by \eqref{KLE}  and  Theorem 2.6, for $z\in D_t$, 
\begin{eqnarray}
 \dot\Phi_t(z) &=& \dot g_t^0(f_t(z))+(g_t^0)'(f_t(z))\dot f_t(z)\nonumber\\
&= &\frac{2\Phi_t'(\xi(t))^2}{g_t^0(f_t(z))-U(t)} +(g_t^0)'(f_t(z))\cdot 2\pi f_t'(z)\Psi_{\s(t)}(z.\xi(t))\nonumber\\
&= &\frac{2\Phi_t'(\xi(t))^2}{\Phi_t(z)-\Phi_t(\xi(t))}+2\pi\Phi_t'(z)\Psi_{\s(t)}(z,\xi(t)). \label{e:13}
\end{eqnarray}

In view  of Lemma 2.5, by an argument similar to that in the paragraphs below (6.32) of \cite{CF}, we can deduce from \eqref{e:13} that $\Phi_t(z)$ is differentiable in $t$ for every  $z\in \partial\HH$, and  $\dot\Phi_t(z)$ is jointly continuous in $(t, z)\in (0, \infty) \times \partial \HH.$   Since $\xi(t)$ is the solution of the SDE \eqref{SDE1},
the above joint continuity together with Lemma \ref{L:2.5} allows us to 
 apply a generalized It\^o formula to $U_t = \Phi_t (\xi (t))$;  see Remark \ref{R:2.8} below.
We thus get 
\begin{eqnarray*}
 dU(t) &=& \dot\Phi_t(\xi(t))dt +\Phi'_t(\xi(t)) \left( \alpha(\s(t)-\wh\xi(t))dB_t+b(\s(t)-\wh\xi(t))dt \right)  \\
 && +\frac12 \Phi_t''(\xi(t))\alpha(\s(t)-\wh\xi(t))^2dt
\end{eqnarray*}
An argument similar to that in the paragraphs below (6.32) of \cite{CF} also yields the identity 
$$
\dot \Phi_t(\xi(t))=\lim_{z\to \xi(t),\; z\in D_t}\dot \Phi_t(z).
$$
Rewriting the right hand side of \eqref{e:13} as
\[
\left( \frac{2\Phi_t'(\xi(t))^2}{\Phi_t(z)-\Phi_t(\xi(t))}-\frac{2\Phi_t'(\xi(t))}{z-\xi(t)}\right)
+2\pi \Phi_t'(\xi(t))\left(\Psi_{\s(t)}(z,\xi(t))+\frac{1}{\pi}\frac{1}{z-\xi(t)}\right),
\]
we obtain from \eqref{e:13} and \cite[Lemma 6.1]{CF}
\[\dot\Phi_t(\xi(t))=-3\Phi_t''(\xi(t))+\Phi_t'(\xi(t))\; b_{\rm BMD}(\s(t)-\wh \xi(t)).\]
Therefore
\begin{eqnarray*}
dU(t) &=& \left( -3\Phi_t''(\xi(t))+\Phi_t'(\xi(t))b_{\rm BMD}(\xi(t)-\wh\xi(t))\right) dt\\
&&+\Phi'_t(\xi(t)) \left( \alpha(\s(t)-\wh\xi(t))dB_t+b(\s(t)-\wh\xi(t))dt \right)
+\frac12 \Phi_t''(\xi(t))\alpha(\s(t)-\wh\xi(t))^2dt ,
\end{eqnarray*}
which is  \eqref{e:12}. \qed

\begin{remark}\label{R:2.8}
\rm {\bf (A generalized It\^o formula)}\quad
Exercise (IV.3.12) in the book \cite{RY} formulates a generalized It\^o formula  
 for $g(X_t,\omega,t)$, 
the composition  of an adapted random function $g(x,\omega,t)$,
$x\in \RR,\, t\ge 0,$ and a continuous semimartingale $X$.    
We like to point out that in addition to the conditions i), ii), iii) and iv)  stated 
in \cite[Exercise IV.3.12]{RY}, the following condition 

\medskip\noindent 
v)\ $g_x(x,\omega,t),\ g_{xx}(x,\omega,t)$ and 
$g_t(x,\omega,t)$ are locally bounded in $(x,t)$ 

\medskip\noindent 
should also be required  
 for the validity of the  generalized It\^o formula 
(a private communication by Masanori Hino).  Of course, 
if these   partial derivatives  are jointly continuous in $(x, t)$, then 
condition v) is satisfied.   
This type of generalized It\^o formula has been frequently utilized in the literatures on SLE by referring to \cite[(IV.3.12)]{RY} but without verifying condition v) which is by no means trivial.    
This is part of  
the reasons 
why we spent considerable efforts in \cite{CF} to establish the joint continuity of  certain  functions
such as those summarized  in Lemma \ref{L:2.5} and 
of the function  
$\dot \Phi_t(z),\; z\in \partial\HH,$  derived from the identity \eqref{e:13}.
\qed
\end{remark}

\section{Proof of property (C)}\label{S:3}

In this section, we present a proof of Proposition 2.4, using the probabilistic representation of $\Im g_t(z)$ in \cite[Theorem 7.2]{CFR} as well as that of $\Im g_t^0(z)$ obtained from \cite[Theorem 7.2]{CFR} by taking  $D=\HH$.

Recall that $g_t(z),\; t\in [0,t_z),$ is the unique solution of \eqref{KLE} 
with the maximal interval $[0,t_z)$ of existence, and 
$F_t=\{z\in D: t_z\le t\}$. We know that 
$g_t(z)$ is continuous in $t$,   and $g_t$ is 
the canonical map from  $D\setminus F_t$. 
 Let $G_t=\{z\in D: t_z<t\}.$  Then
\begin{equation}\label{eq:9}
\bigcap_{s>t} F_{s}=F_t,\qquad \bigcup_{s<t} F_{s}=G_t.
\end{equation}

Let $g_t^0$ be 
the canonical Riemann map from $\HH\setminus F_t$.
By virtue of Theorem 7.2 of \cite{CFR} with  $D=\HH$ (see also \cite[(3.5)]{L1}), $\Im g_t^0(z)$ admits the expression
\begin{equation}\label{eq:10}
\Im g_t^0(z)=\Im z - \E_z^\HH\left[\Im Z_{\sigma_{F_t}}^\HH: \sigma_{F_t}<\infty\right],\ z\in \HH\setminus F_t,
\end{equation}
where $ Z^\HH=(Z_t^\HH, \zeta^\HH, \P_z^\HH)$ is the absorbing Brownian motion (ABM) on $\HH$, 
and $\sigma_{F_t}:=\inf\{s>0: Z^\HH_s \in F_t\}$.

\begin{lem} \label{L:3.1} \ Fix $a>0.$  $\Im g_t^0(z)$ is continuous in $t\in [0,a]$ for each $z\in \HH\setminus F_a$ if and only if
\begin{equation}\label{eq:11}
\E_z^\HH\left[\Im Z_{\sigma_{G_t}}^\HH; \sigma_{G_t}<\infty\right]
=\E_z^\HH\left[\Im Z_{\sigma_{F_t}}^\HH; \sigma_{F_t}<\infty\right],
\end{equation} 
for $t\in (0, a]$ and $ z\in \HH\setminus F_a.$
\end{lem}

\pf\ Since $\sigma_{F_{s}} \downarrow \sigma_{G_t}$ as $s \uparrow t$ by 
\eqref{eq:9} (cf. \cite[Chapter 1, (10.4)]{BG}), 
we see from \eqref{eq:10} that \eqref{eq:11} is equivalent to the left continuoity of $\Im g_t^0(z)$ in $t$.
On the other hand, $\Im g_t^0(z)$ is right continuous in $t$
 because $g_t^0(z)$ is right differentiable in $t$ by Proposition \ref{P:2.3}.  \qed

Let $K=\bigcup_{j=1}^N C_j$ and $v_t^*(z)=\Im g_t(z)$.  
Denote by $Z^{\HH,*}=(Z^{\HH,*}_t, \P^{\HH,*}_z)$
the BMD on $D^*=D\cup K^*$ with $K^*=\{c_1^*,\dots, c_N^*\}$ obtained from 
the ABM $Z^\HH$ by shorting each slit $C_i$ as a single point $c_i^*.$
According to \cite[Theorem 7.2]{CFR}, $v_t^*(z)$ can be expressed in terms of the ABM $Z^\HH$ and BMD $Z^{\HH,*}$ as follows:
\begin{equation}\label{eq:12}
v_t^*(z)=v_t(z)+\sum_{j=1}^N \P_z^\HH \left(\sigma_K<\sigma_{F_t},\; Z_{\sigma_K}^\HH\in C_j \right)
v_t^*(c_j^*),
\quad  z\in D\setminus F_t,
\end{equation}
where
\begin{equation}\label{eq:13}
v_t(z)=\Im z-\E_z^\HH\left[ \Im Z_{\sigma_{F_t\cup K}}^\HH; \sigma_{F_t\cup K}<\infty \right],
\end{equation}
\begin{equation}\label{e:7.4}
v_t^*(c_i^*)=\sum_{j=1}^N\frac{M_{ij}(t)}{1-R_{i}^*(t)}\;
\int_{\eta_j}v_t(z)\nu_j (dz),
\quad  1\le i\le N.
\end{equation}
Here $\eta_1,\cdots, \eta_N$ are mutually disjoint smooth Jordan curve
 surrounding $C_1,\cdots, C_N$, respectively,
\begin{equation}\label{e:7.5}
\nu_i(dz)=\P_{c_i^*}^{\HH,*}\left(Z_{\sigma_{\eta_i}}^{\HH,*}\in dz\right),
\quad 1\le i\le N,
\end{equation}
\begin{equation}\label{e:7.6}
\ R_{i}^*(t)=\int_{\eta_i}\P_z^\HH\left(\sigma_K<\sigma_{F_t},
\ Z_{\sigma_K}^\HH\in C_i\right)\nu_i(dz),\quad 1\le i\le N,
\end{equation}
and $M_{ij}(t)$ is the $(i,j)$-entry of the matrix
$M(t)=\sum_{n=0}^\infty (Q^*(t))^n$ for a matrix $Q^*(t)$ with
entries
\begin{equation}\label{e:7.7}
q_{ij}^*(t) =
\begin{cases}
\P_{c_i^*}^{\HH,*} \big(\sigma_{K^*}<\sigma_{F_t},\; Z_{\sigma_{K^*}}^{\HH,*}=c_j^*\big) / (1-R_{i}^*(t))
&\hbox{\rm if  } i\neq j,\\
 0 &\hbox{\rm if } i=j,
\end{cases}
\quad  1\le i,j\le N.
\end{equation}

\begin{lem} \label{L:3.2} \  For every $1\leq j\leq N$, $ v_t^*(c_j^*) >0$ for every $t>0$ and 
\begin{equation}\label{eq:14}
 \sup_{0\le t\le a} v_t^*(c_j^*)<\infty \quad \text{\rm for each}\ a>0  .
\end{equation}
\begin{equation}\label{eq:14b}
  v_t^*(c_j^*) >0,\quad t>0,\quad 1\le j\le N.
\end{equation}
\end{lem}

\pf\ For $0\le t\le a $  and $1\le i\le N$, let
\[\lambda_i(t)= \sum_{j=1}^N q_{ij}^*(t) 
\quad \hbox{and} \quad 
 \gamma_i(t)=\P_{c_i^*}^{\HH,*}\left( \sigma_{K^*}<\sigma_{F_t},\; Z_{\sigma_{K^*}}^{\HH,*} \neq c_i^*\right), \quad 1\le i\le N.
\]
Note that  $\lambda_i(t)=\gamma_i(t)/(1-R_i^*(t))$ and
\begin{equation}\label{eq:14c} 
1-R_i^*(t)=\gamma_i(t)+
\int_{\eta_i} P_z^\HH\left(\sigma_{F_t}<\sigma_K\right) \nu_i(dz) 
+\int_{\eta_i} \P_z^\HH\left(\sigma_{F_t\cup K}=\infty \right) \nu_i(dz) .
\end{equation}
Therefore 
\begin{eqnarray*}
1-\lambda_i(t) &=& \frac{1-R^*_i(t)-\gamma_i(t)}{1-R_i^*(t) }
\geq \int_{\eta_i}  \P_z^\HH(\sigma_{F_t\cup K}=\infty) \nu_i(dz) \\
&\geq & \inf_{1\leq j\leq N}  \int_{\eta_j}  \P_z^\HH(\sigma_{F_t\cup K}=\infty) \nu_j (dz) =:\delta_0>0 .
\end{eqnarray*}
Hence $\lambda_i(t)\le 1-{{\delta_0}}$. Consequently,  
$(Q^*(t))^n{\bf 1}\le (1-{{\delta_0}})^n{\bf 1}$ and so $ M{\bf 1}\le  \delta_0^{-1} {\bf1}$.
  Therefore we have by \eqref{e:7.4} and \eqref{eq:14c} that 
$ v_t^*(c_i^*)\le  \sum_{j=1}^N  \delta_0^{-2}  m_j$ for all $t\in [0, a]$, 
where $m_j$ is the maximum of the $y$-th coordinate of points in $\eta_j.$

On the other hand, \eqref{e:7.4} implies $v_t^*(c_i^*)\ge \int_{\eta_i}v_t(z)\nu_i(dz).$  In view of \eqref{eq:13}, $v_t(z)$ is a non-negative harmonic function on $\HH\setminus (F_t\cup K)$ that is strictly positive when $\Im z$ is large.
Hence $v_t(z)>0$ for any $z\in \HH\setminus (F_t\cup K)$ and $t>0$, 
yielding \eqref{eq:14b}.   \qed

\begin{prop} \label{P:3.3} \ The identity {\rm \eqref{eq:11}} holds, and so $g^0_t (z) $ is continuous in $t\in [0, a]$
for every $z\in \HH \setminus \FF_a$ and $a>0$.
\end{prop}

\pf\ Note that $v_t^*(z)=\Im g_t(z)$ is continuous in $t$ since so is  $g_t(z)$.
By \eqref{eq:12}-\eqref{eq:13}, for $z\in D\setminus F_t$, 
\begin{equation} \label{e:3.12}
v_t^*(z)  = \Im z - \E_z^\HH\left[ \Im Z_{\sigma_{F_t\cup K}}^\HH; \sigma_{F_t\cup K}<\infty\right]
 +\sum_{j=1}^N \P_z^\HH \left( \sigma_K<\sigma_{F_t}, Z_{\sigma_K}^\HH\in C_j \right) v_t^*(c_j^*) .
 \end{equation}
For each fixed $t\in (0,a]$ and any sequence $t_n$ increasing to $t$, by \eqref{eq:14},
there is a subsequence $t_{n_k}$ such that $\lim_{k\to\infty} v_{t_{n_k}}^*(c_j^*)=a_j \in [0, \infty)$.  
Since $ F_{t_n} \uparrow  G_t $, we have
\begin{equation} \label{e:3.13}
v_t^*(z)=\lim_{k\to\infty} v_{t_{n_k}}^* (z)
= \Im z - \E_z^\HH\left[ \Im Z_{\sigma_{G_t\cup K}}^\HH; \sigma_{G_t\cup K}<\infty\right]
 +\sum_{j=1}^N \P_z^\HH \left( \sigma_K<\sigma_{G_t}, Z_{\sigma_K}^\HH\in C_j \right)   a_j .
\end{equation}
 Taking $z\to C_j$ in \eqref{e:3.12} and \eqref{e:3.13} yields   $a_j=v_t^*(c_j^*)$ for each $1\leq j\leq N$. 
Thus we have from \eqref{e:3.12} and \eqref{e:3.13} that
\begin{eqnarray}
&& \E_z^\HH\left[ \Im Z_{\sigma_{G_t\cup K}}^\HH; \sigma_{G_t\cup K}<\infty\right]
-   \E_z^\HH\left[ \Im Z_{\sigma_{F_t\cup K}}^\HH; \sigma_{F_t\cup K}<\infty\right]  \nonumber \\
&=&\sum_{j=1}^N  \left( \P_z^\HH \left(\sigma_K<\sigma_{G_t}, Z_{\sigma_K}^\HH\in C_j \right) 
 - \P_z^\HH \left(\sigma_K<\sigma_{F_t}, Z_{\sigma_K}^\HH\in C_j \right) \right) v_t^*(c_j^*) .     \label{e:3.14} 
 \end{eqnarray} 

Each term on the right hand side of \eqref{e:3.14} is non-negative since $G_t =F_{t-}\subset F_t$. 
 On the other hand, $\Im z$ is a positive harmonic in $\HH$
and so $\Im Z^\HH_t$ is a non-negative supermartingale.
By the optional  sampling theorem, we have for every $z\in \HH$
and any stopping time $T$, we have
\begin{equation}\label{e:3.16} 
 \Im z \geq \E_z^\HH  [ \Im Z^\HH_T; T <\infty].
\end{equation}
Since $\sigma_{G_t\cup K} \geq \sigma_{F_t\cup K}$, 
we have
$$ \Im z  \geq \E_z^\HH\left[ \Im Z_{\sigma_{F_t\cup K}}^\HH; \sigma_{F_t\cup K}<\infty\right]
\geq \E_z^\HH\left[ \Im Z_{\sigma_{G_t\cup K}}^\HH; \sigma_{G_t\cup K}<\infty\right] \geq 0  , 
$$
where in the second inequality we used the strong Markov property
of $Z^\HH$ at stopping time $\sigma_{F_t\cup K}$ and \eqref{e:3.16}.
Thus both sides of \eqref{e:3.14} have to be identically zero. 
 As $v_t^*(c_j^*)>0$ for each $1\le j\le N$   
by \eqref{eq:14b}, 
we must   have for $z\in D\setminus F_t$, 
\begin{equation}\label{eq:18}
\E_z^\HH\left[\Im Z_{\sigma_{F_t\cup K}}^\HH; \sigma_{F_t\cup K}<\infty\right]
=\E_z^\HH\left[\Im Z_{\sigma_{G_t\cup K}}^\HH; \sigma_{G_t\cup K}<\infty\right],
\end{equation}
and
\begin{equation}\label{eq:19}
\P_z^\HH \left(\sigma_K<\sigma_{G_t}, Z_{\sigma_K}^\HH\in C_j \right) 
 = \P_z^\HH \left(\sigma_K<\sigma_{F_t}, Z_{\sigma_K}^\HH\in C_j \right)  
 \quad \hbox{for every } 1\leq j\leq N.
\end{equation}
It follows from the above two displays  that for $z\in \HH\setminus (K\cup F_t)$, 
\begin{equation}\label{eq:20}
  \P_z^\HH \left(\sigma_K<\sigma_{F_t} \right) = \P_z^\HH \left(\sigma_K<\sigma_{G_t} \right)  
\quad \hbox{and} \quad 
\E_z^\HH\left[\Im Z_{\sigma_{F_t}}^\HH; \sigma_{F_t}<\sigma_K\right]=\E_z^\HH\left[\Im Z_{\sigma_{G_t}}^\HH; \sigma_{G_t}<\sigma_K\right] .
\end{equation}

Take a bounded smooth domain $V\subset \HH$ such that $K\subset V$ and $V\cap F_t=\emptyset$.
Let $\Gamma=\partial V.$   Define
$\sigma_1=\sigma_K$, $ \tau_1 = \inf\{t\geq \sigma_1: Z_t^\HH \in \Gamma\} $,
 and for $n\geq 1$,
$$ \sigma_{n+1} = \inf\{ t> \tau_n: Z_t^\HH \in K \} , \quad 
      \tau_{n+1} = \inf\{ t> \sigma_{n+1}: Z_t^\HH \in \Gamma \}.
 $$
We claim that the following holds for every $n\geq 1$ and  $z\in \HH\setminus (K\cup F_t)$, 
\begin{equation}\label{eq:21}
\P_z^\HH(\sigma_n<\sigma_{F_t})=\P_z^\HH(\sigma_n<\sigma_{G_t})
 \quad \hbox{and} \quad 
\P_z^\HH(\tau_n<\sigma_{F_t})=\P_z^\HH(\tau_n<\sigma_{G_t}) .
\end{equation}
We prove this by induction.  Clearly the first identity in \eqref{eq:21} holds for $n=1$ by   \eqref{eq:20},
while by the continuity  of the sample paths of $Z^\HH$, 
$$ 
\P_z^\HH (\tau_1<\sigma_{F_t} ) =\P_z^\HH ( \sigma_K<\sigma_{F_t} )= \P_z^\HH ( \sigma_K<\sigma_{G_t} )
=\P_z^\HH(\tau_1<\sigma_{G_t}).
$$
So \eqref{eq:21} holds for $n=1$. 
 Assume that \eqref{eq:21} holds for $n\geq 1$.
  Then   by the strong Markov property  of $Z^\HH$ and \eqref{eq:20}, 
\begin{eqnarray*}
&&\P_z^\HH(\sigma_{n+1}<\sigma_{F_t})=\P_z^\HH \left(\tau_n+\tau_K \circ\theta_{\tau_n}<\sigma_{F_t}  , \  \tau_n<  \sigma_{F_t}\right)\\
& =& \E_z^\HH\left[\P_{Z_{\tau_n}}^\HH(\sigma_K <\sigma_{F_t}); \tau_n<\sigma_{F_t}\right]
=\E_z^\HH\left[\P_{Z_{\tau_n}}^\HH(\sigma_K<\sigma_{G_t}); \tau_n<\sigma_{G_t}\right]
=\P_z^\HH(\tau_{n+1}<\sigma_{G_t} ),
\end{eqnarray*}
and by the continuity of $Z^\HH$, 
\begin{eqnarray*}
&&\P_z^\HH(\tau_{n+1}<\sigma_{F_t})=\P_z^\HH \left(\sigma_{n+1}+\tau_{\Gamma} \circ \theta_{\sigma_{n+1}}<\sigma_{F_t} \right)
=  \P_z^\HH \left(\sigma_{n+1} <\sigma_{F_t} \right) \\
& =& \P_z^\HH \left(\sigma_{n+1} <\sigma_{G_t} \right)
=\P_z^\HH(\tau_{n+1}<\sigma_{G_t} ).
\end{eqnarray*}
Hence \eqref{eq:21} holds for $n+1$ and so for all $n\geq 1$ by induction.

Now, by the strong Markov property of $Z^\HH$, \eqref{eq:21} and \eqref{eq:20},
 we have for $z\in \HH\setminus (K\cup F_t)$
\begin{eqnarray*}
 \E_z^\HH\left[ \Im Z_{\sigma_{F_t}}^\HH; \sigma_{F_t}<\infty\right] 
& =& \E_z^\HH\left[ \Im Z_{\sigma_{F_t}}^\HH; \sigma_{F_t} <\sigma_K\right]
+\sum_{n=1}^\infty\E_z^\HH\left[ \Im Z_{\sigma_{F_t}}^\HH; \sigma _n<\sigma_{F_t} <\sigma_{n+1}  \right]\\
& = & \E_z^\HH\left[ \Im Z_{\sigma_{F_t}}^\HH; \sigma_{F_t}<\sigma_K\right]
+\sum_{n=1}^\infty\E_z^\HH\left[ \E_{Z_{\tau_n}}^\HH[\Im Z_{\sigma_{F_t} }^\HH; \sigma_{F_t} <\sigma_K] ; \tau_n<\sigma_{F_t} \right]
\\
 & = & \E_z^\HH\left[ \Im Z_{\sigma_{G_t}}^\HH; \sigma_{G_t}<\sigma_K\right]
+\sum_{n=1}^\infty\E_z^\HH\left[ \E_{Z_{\tau_n}}^\HH[\Im Z_{\sigma_{G_t} }^\HH; \sigma_{G_t} <\sigma_K] ; \tau_n<\sigma_{G_t} \right]
\\
& =& \E_z^\HH\left[ \Im Z_{\sigma_{G_t}}^\HH; \sigma_{G_t} <\sigma_K\right]
+\sum_{n=1}^\infty\E_z^\HH\left[ \Im Z_{\sigma_{F_t}}^\HH; \sigma _n<\sigma_{G_t} <\sigma_{n+1}  \right] \\ 
&=& \E_z^\HH\left[ \Im Z_{\sigma_{G_t}}^\HH; \sigma_{G_t}<\infty\right] . 
\end{eqnarray*}
This establishes \eqref{eq:11}. The rest of the claim follows from Lemma \ref{L:3.1}.  \qed

For $0\le s<t\le a,$  define  $g_{t,s}^0=g_s^0\circ(g_t^0)^{-1}$,
 which is a conformal map from $\HH$ onto $\HH\setminus g_s^0(F_t\setminus F_s)$. Its inverse 
$(g_{t,s}^0)^{-1}$ is 
the canonical Riemann map from $\HH\setminus g_s^0(F_t\setminus F_s)$.
 Let $\ell_{t,s}$ be the set of all limitting points of $(g_{t,s}^0)^{-1}\circ g_s^0(z)=g_t^0(z)$ as $z$ approaches to $F_t\setminus F_s.$   Then $\ell_{t,s}$ is a compact subset of $\partial\HH$ and $(g_{t,s}^0)^{-1}$ sends $\partial\HH\setminus \overline{g_s^0(F_t\setminus F_s)}$ into $\partial\HH$ homeomorphically. 

Let $\Lambda=\{x+iy:a<x<b,\;0<y<c\}$ be a finite rectangle such that $\ell_{t,s} \subset \{x+i0+: a<x<b\}.$   Then $\Im g_{t,s}^0(z)\le \Im (g_t^0)^{-1}(z)$ by \eqref{eq:10} that is uniformly bounded in $z\in \Lambda$ so that it admits finite limit
\begin{equation}
\Im g_{t,s}^0(x+i0+)=\lim_{y\downarrow 0} \Im g_{t,s}^0(x+iy)\quad \text{\rm for a.e.}\ x\in (a,b).
\end{equation}

The following lemma can be established in a similar way as that of \cite[Lemma 6.3]{CF}. We omit its proof here.

\begin{lem} \label{L:3.4} \ For $0\le s<t\le a,$ it holds that
\begin{equation}
a(t)-a(s)=\frac{1}{\pi}\int_{\ell_{t,s}} \Im g_{t,s}^0(+i0+)~dx, 
\end{equation}
\begin{equation}
g_t^0(z)-g_s^0(z)=-\frac{1}{\pi}\int_{\ell_{t,s}}\frac{1}{g_t^0(z)-x} \Im g_{t,s}^0(x+i0+)~dx,\quad z\in \HH\setminus F_t.
\end{equation}
\end{lem}

\medskip

\noindent
{\bf Proof of Proposition 2.4. } \ We know from Proposition \ref{P:3.3} that  $\Im g_t^0(z)$ is continuous in $t\in [0,a]$ for each $z\in \HH\setminus F_a$.  As $\Im g_t^0(z)$ is harmonic in $z\in \HH\setminus F_a$, it is jointly continuous in $(t,z)\in [0,a]\times (\HH\setminus F_a).$
By Lemma 3.3, we have  
\[|g_s^0(z)|\le |g_t^0(z)|+ \sup_{x\in \ell_{t,0}}\frac{a_t}{|g_t^0(z)-x|},\quad s\in [0,t].\]
Therefore we can show as in the proof of \cite[Theorem 7.4]{CFR} that $g_t^0(z)$ is locally equi-continuous and locally uniformly bounded.  The joint continuity of $g_t^0(z)$ then follows as in the proof of \cite[Lemma 6.5]{CF}. \qed

\section{Basic relations between ${\rm SKLE}_{\alpha,b}$ and  SLE}
 
In view of \cite[(7.20)]{CF} 
applied to  the case $D=\HH$ (see also \cite[(3.7)]{L1}), the half-plane capacity $a(t)$ of the hull $F_t$ relative to $g_t^0$ admits the expression
\[a(t)=\frac{2R}{\pi}\int_0^\pi \E_{R e^{i\theta}}^\HH\left[\Im Z_{\sigma_{F_t}}^\HH: \sigma_{F_t}<\infty\right]d\theta,\]
in terms of the ABM $Z^\HH$ on $\HH$ for $R>0$ with $F_t\subset \{z\in \HH:|z|<R\}.$  Since the SKLE $\{F_t\}$ is strictly increasing in $t$ by virtue of \cite[Theorem 5.8]{CF}, we can see as in the proof of \cite[Lemma 5.15]{CF} that $a(t)$ is strictly increasing in $t$.

By Lemma 2.2 and Lemma 2.5,
\begin{equation}\label{3.1}
a(t)=2\int_0^t |\Phi_s'(\xi(s))|^2 ds.
\end{equation}  
We  reparametrize the SKLE hulls $\{F_t\}$ by the inverse function $a^{-1}$ of $a$ and define
\begin{equation}\label{3.2}
\c F_t=F_{a^{-1}(2t)},\qquad 0\le t < \tau_0:=a(\infty)/2.
\end{equation}
Accordingly, the associated Riemann maps $\{g_t^0\}$ and the process $U(t)$ are time changed into
\begin{equation}\label{3.3}
\c g_t^0=g^0_{a^{-1}(2t)}, \quad \c U(t)=U(a^{-1}(2t)),\quad 0\le t< \tau_0 .
\end{equation}
It then follows from \eqref{e:11} that $z(t)=\c g^0_t(z)$ is a solution of the Loewner equation
\begin{equation}\label{3.4}
\frac{d}{dt} z(t)\;=\;\frac{2}{z(t)-\c U(t)},\quad z(0)=z \in \HH.
\end{equation} 

\begin{thm}\label{T:4.1} 
 $\{\c F_t; t\in [0, \tau_0) \}$ has the same law as the Loewner evolution driven by the path of the continuous process $\c U(t)$
 up to the random time $\tau_0$;   namely, for the unique solution $z(t),\; 0\le t<t_z,$ of {\rm \eqref{3.4}},
\begin{equation}\label{3.5}
\{\c F_t;  t\in [0, \tau_0)\}   \hbox{ has the same distribution as }   \{\{z\in \HH: t_z\le t\}; t\in [0, \tau_0) \}.
\end{equation}
\end{thm}

Let $M_t=\int_0^t \Phi_s'(\xi(s))dB_s$.  By \eqref{3.1}, $\< M\>_t=\int_0^t\Phi_s'(\xi(s))^2ds = a(t)/2$ so that $\c B_t=M_{a^{-1}(2t)}$ is a Browninan motion.
The formula \eqref{e:12} can be rewritten as
\begin{eqnarray}\label{3.6}
 \check U (t)&=&
 \xi(0) + \int_0^t\wt \Phi_s'(\wt\xi(s))^{-1}\left( b(\wt\s(s)-\wh{\wt\xi}(s))+b_{\rm BMD}(\wt\s(s)-\wh{\wt\xi}(s))\right) ds\nonumber\\
&&+\frac12 \int_0^t \wt \Phi''_s(\wt\xi(s))\cdot \wt \Phi_s'(\wt\xi(s))^{-2}  \left( \alpha(\wt\s(s)-\wh{\wt\xi}(s))^2-6 \right) ds\nonumber\\
&& 
+ \int_0^t \alpha(\wt \s(s)-\wh{\wt\xi}(s)) d\check B_s,
\end{eqnarray}
where
$\wt \Phi_s'(z):=\Phi_{a^{-1}(2s)}'(z)$, $\Phi_s''(z):=\Phi_{ a^{-1}(2s)}''(z)$,
$\wt\xi(t ):= \xi( a^{-1}(2t))$ and $\wt\s_j( t)= \s_j( a^{-1}(2t))$
for  $1\le j\le 3N$.
Note that since $\Phi_t(z)$ is univalent in $z$ on the region $D_t\cup \Pi D_t\cup \partial \HH$,  $\Phi_t'(z)$ never vanishes there.
\eqref{3.6} particularly means that $\c U(t)$ is a continuous semimartingale. 

From Theorem \ref{T:4.1}  and the identity \eqref{3.6}, we can obtain immediately the following two theorems.

\begin{thm}\label{T:4.2}  ${\rm SKLE}_{\sqrt{6}, -b_{\rm BMD}}$ being reparametrized as {\rm \eqref{3.2}} has the same distribution as ${\rm SLE}_6$ over the time  interval $[0, \tau_0 ).$
\end{thm}

\begin{thm}\label{T:4.3}  For a positive constant $\alpha,$ there exists a sequence of hitting times $\{\sigma_n\}$ increasing to $\tau_0$ such that ${\rm SKLE}_{\alpha, b}$ being reparametrized as {\rm \eqref{3.2}} has the same distribution as ${\rm SLE}_{\alpha^2}$ over each time  interval $[0,\sigma_n]$ under a suitable Girsanov transform.
\end{thm}

When $\alpha$ is a positive constant, it follows from Theorem \ref{T:4.3} and \cite{RS} that ${\rm SKLE}_{\alpha,b}$ is generated by a continuous curve $\gamma$ and that $\gamma$ is simple when $\alpha\le 2$, self-intersecting when $2<\alpha \leq 2\sqrt{2}$ and space-filling
when $\alpha > 2\sqrt{2}$.  

  \section{Locality property of ${\rm SLE}_6$ in several canonical domains}\label{S:5}

It has been demonstrated in \cite[Theorem 6.11]{CF}
that ${\rm SKLE}_{\alpha, -b_{\rm BMD}}$ 
enjoys the locality property for a positive constant $\alpha$ if and only if 
$\alpha=\sqrt{6}.$   The proof is being carried out independently of the locality of ${\rm SLE}_6$.
The next subsection will concern the question:

\smallskip\noindent
{\bf (Q)} \ Is there any alternative proof of the locality of 
 ${\rm SKLE}_{\sqrt{6}, -b_{\rm BMD}}$ based on Theorem 4.2 ?

\subsection{Locality of chordal ${\rm SLE}_6$ and ${\rm SKLE}_{\sqrt{6},-b_{\rm BMD}}$}

Let $\Phi$ be a locally real conformal transformation from an $\HH$-neighborhood $\cal N$  of a subset of $\partial\HH$ into $\HH$ in the sense of \cite[\S 4.6]{L1}.
Theorem 6.13 of \cite{L1} claimed a {\it locality} of ${\rm SLE}_6$ {\it relative to} $\Phi$ in the following sense: the ${\rm SLE}_6$-hulls $\{K_t\}$ have the same law as $\{\Phi(K_t)\}$ until the exit time from $\Phi ({\cal N} \cup \partial \HH)$ up to a time change.   
The proof was based on a generalized Loewner equaiton
\begin{equation}\label{4.7}
\frac{dg_t^*(z)}{dt} = \frac{2 \Phi_t'(\xi(t))^2}{g_t^*(z)- U^*(t)},\quad g_0^*(z)=z,\quad U^*(t)=\Phi_t(\xi(t)),
\end{equation}
for the canonical Riemann map $g^*_t(z)$ from $\HH\setminus \Phi(K_t).$
 Here $d\xi(t)=\sqrt{6}dB_t$ and 
\begin{equation}\label{4.8}
\Phi_t\, :=\; g_t^*\circ \Phi\circ g_t^{-1},
\end{equation} 
where $g_t(z)$ is  the solution of the Loewner equation \eqref{1.1}.

But the equation \eqref{4.7} was rigorously proved in \cite{L1} only in the right derivative sense just as the proof of Proposition \ref{P:2.3} of this paper.  In order to make it a genuine ODE, we need to verify the joint continuity of $g^*_t(z)$ in $(t,z),$ which can be shown when $\Phi$ is 
the canonical Riemann map $\varphi_A$ from $\HH\setminus A$ for any 
$\HH$-hull $A\subset \HH$ by using 
the probabilistic representation of $\Im \Phi_t (z)$. 
Indeed, in this case, we have 
\begin{equation}\label{4.9}
\Phi_t(z)= \varphi_{g_t(A)}(z),\quad z\in \HH\setminus g_t(A),
\end{equation}
for the canonical Riemann map $\varphi_{g_t(A)}$ from $\HH\setminus g_t(A)$
and so  $\Im \Phi_t(z)$ admits a probabilistic expression
\begin{equation}\label{4.10}
\Im \Phi_t(z)=\Im z-\E_z^\HH\left[ \Im Z_{\sigma_{g_t(A)}}^\HH; \sigma_{g_t(A)}<\infty\right]
\end{equation}
in terms of the ABM $(Z_t^\H, \P_z^\HH)$ on $\HH$ in view of \cite[(3.5)]{L1}.  Define
\begin{equation}\label{4.11}
q_t(z) = \Im g_t(z)-\E_z^\HH\left[ \Im g_t(Z_{\sigma_A}^\HH); \sigma_A<\infty\right],\quad z\in \HH\setminus (F_t\cup A) .
\end{equation}
Due to the invariance of the ABM under the conformal map $g_t$, we have $\Im \Phi_t(g_t(z))= q_t(z).$  Since $g^*_t=\Phi_t\circ g_t\circ \varphi_A^{-1}$ by \eqref{4.8}, we obtain for each $T <\tau_A:= \inf\{t: \overline K_t\cap \overline A\neq \emptyset\}$  
\begin{equation}\label{4.12}
\Im g^*_t(z)= q_t(\varphi_A^{-1}(z)),\quad t\in [0,T],\quad z\in \HH\setminus 
\varphi_A(K_T).
\end{equation}
As $g_t(z)$ is the solution of the Loewner equation \eqref{1.1}, $\Im g_t(z)$ is jointly continuous and bounded by $\Im z$.  Hence $q_t(z)$ is continuous in $t$ for each $z\in \HH\setminus K_t\setminus A$ by \eqref{4.11} and so is $\Im g^*_t(z)$ for each $z\in \HH\setminus \varphi_A(K_T).$   This continuity implies the joint continuity of $g^*_t(z)$ just as in the last part of Section 3 and so \eqref{4.7} becomes a genuine ODE.  Using a generalized It\^o formula as the proof of Theorem \ref{T:2.7}, we can then obtain $dU^*(t)=\sqrt{6}\Phi_t'(\xi(t))dB_t,\ \  t<T$.

 We have thus given a first rigorous proof of 
the locality of ${\rm SLE}_6$ relative to $\varphi_A$.

\begin{prop}{\bf (Locality of chordal ${\rm SLE}_6$ in the sense of \cite{LSW3}).}\quad For any $\HH$-hull $A$, let $\varphi_A$ be the canonical Riemann map from $\HH\setminus A.$  Then the image hulls $\{\varphi_A(K_t)\}$ of ${\rm SLE}_6$-hulls $\{K_t\}$ has under a reparametrization the same law as $\{K_t\}$ up to the first hitting time $\sigma_A$ of $A$.  
\end{prop}

Notice that, in view of \cite{RS}, ${\rm SLE}_6$-hulls $\{K_t\}$ are generated by continuous self-interesecting curves, thus so are the image hulls $\{\varphi_A(K_t)\}.$  Accordingly, the classical argument for a Jordan arc yielding the left continuity in $t$ of $g_t^*(z)$ (see \cite[\S 6]{CFR}) cannot be applied and 
no proof of the continuity of $g_t^*$ in $t$ seems to be available other than the probabilistic method we employed above. 

Now, for any standard slit domain $D$ and any $\HH$-hull $A\subset D,$ consider the canonical conformal map $\Phi_A$ from $D\setminus A.$  
Note that $\Phi_A$ 
is a specific locally real conformal map from the $\HH$-neighborhood $D\setminus A$ of $\partial \HH\setminus \overline A.$  If we could verify the locality of ${\rm SLE}_6$ relative to $\Phi_A$, then the locality of ${\rm SKLE}_{\sqrt{6},-b_{\rm BMD}}$ could be readily deduced from Theorem 4.2. But $\Phi_t$ defined by \eqref{4.8} for $\Phi=\Phi_A$ does not satisfy \eqref{4.9} unless $D=\HH$ so that the above probabilistic method does not work 
for proving the locality of ${\rm SLE}_6$ relative to $\Phi_A.$ 

So the answer to  question {\bf (Q)} remains negative at present. 
However it may be still possible 
to show the  locality of ${\rm SLE}_6$ relative to $\Phi_A$,
for example, from the point of view that ${\rm SLE}_6$ is the scaling limit of 
the critical percolation exploration process on
triangular lattices,   
although we feel that its rigorous proof would get lengthy.

\subsection{A locality of radial ${\rm SLE}_6$ relative to modified canonical maps}

So far, only chordal SLEs and chordal SKLEs have been considered.  

Consider a linear transformation $\dis\psi(z)=i\frac{1+z}{1-z}$ from the unit disk $\DD=\{z\in \CC:|z|<1\}$ onto $\HH,$ that sends $0$ to $i$ and $1$ to $\infty.$  Its inverse $\psi^{-1}$ sends $\partial\HH$ onto $\partial\DD \setminus \{0\}.$  
Let $\{K_t\}$ and $\{\wh K_t\}$ be the radial ${\rm SLE}_\kappa$ on $\DD$ and the chordal ${\rm SLE}_\kappa$ on $\HH$, respectively.  
A Basic relation of their distributions was investigated 
in \cite[\S 4.2]{LSW2} by using the map $\psi.$

More specifically, define $\wt K_t=\psi^{-1}(\wh K_t).$  $\{\wt K_t\}$ is then a family of random growing hulls on $\DD$ starting at some point of $\partial \DD\setminus \{0\}.$.  In \S 4.2 of \cite{LSW2}, the following statement was established by a right application of a generalized It\^o formula mentioned in Remark 2.9 ($e_t=g_t(1)$ in its proof is a random variable).     

\begin{prop}{\bf (Theorem 4.1 of \cite{LSW2})}\quad The radial ${\rm SLE}_6$ $\{K_t\}$ has under a reparametrization the same law as the $\psi^{-1}$-image $\{\wt K_t\}$ of the chordal ${\rm SLE}_6$ $\{\wh K_t\}$ up to certain hitting time. 
\end{prop} 

For a hull $A$ on $\DD$ with $0\notin A,$ the unique Riemann map $\Phi_A$ from $\DD\setminus A$ onto $\DD$ satisfying $\Phi_A(0)=0,\ \Phi_A'(0)>0,$ is called the {\it canonical map} from $\DD\setminus A.$  We also define a {\it modified canonical map} $\wt \Phi_A$ from $\DD\setminus A$ (onto $\DD$) by
\[\wt \Phi_A=\psi^{-1}\circ \varphi_{\psi(A)}\circ \psi,\]
where $\varphi_{\psi(A)}$ is the canonical Riemann map from $\HH\setminus \psi(A)$ (onto $\HH$).  A modified canonical map $\wt \Phi_A$ is different from the canonical map $\Phi_A$ because $\wt \Phi_A(0)\neq 0$ in general.      
       
The radial ${\rm SLE}_\kappa$ $\{K_t\}$ is said to enjoy the {\it locality property} if, for any hull $A\subset \DD$ with $0\notin A,$ 
$\{\Phi_A(K_t)\}$ has under a reparametrization the same law as $\{K_t\}$ until the hitting time $\tau_A=\inf\{t: \overline K_t\cap \overline A\neq \emptyset\}$ for the canonical map $\Phi_A$ from $\DD\setminus A.$
It readily follows from Proposition 5.1 and Proposition 5.2 that the radial ${\rm SLE}_6$ $\{K_t\}$ enjoys the locality but relative to the modified canonical map $\wt \Phi_A$:

\begin{cor}
For any hull $A\subset \DD$ with $0\notin A,$ 
$\{\wt\Phi_A(K_t)\}$ has under a reparametrization the same law as $\{K_t\}$ until a hitting time not greater than $\tau_A$ for the modified canonical map $\wt\Phi_A$ from $\DD\setminus A.$
\end{cor}

In order to show the locality of the radial ${\rm SLE}_6$ (relative to canonical maps), one may need to make analogous considerations to the proof of Proposition 5.1 first by deriving a generalized Loewner equation in the right derivative sense and then using the absorbing Brownian motion on $\DD.$  We leave its proof to interested readers.  
 
An annulus ${\rm SLE}_\kappa$ was introduced by \cite{Z1} where it was claimed that, for any compact set $F\subset \DD$ containing the origin, the radial ${\rm SLE}_6$ being stopped upon hitting $F$ has the same law as the annulus ${\rm SLE}_6$ up to a time change.  
See the next section.
 A locality of the annulus ${\rm SLE}_6$ can be readily deduced from 
 this special property combined with Corollary 5.3.

\section{K-L equations and SKLEs for other canonical domains}
In this section, we recall and examine Komatu-Loewner equations and stochastic Komatu-Loewner evolutions 
studied in literature for other canonical multiply connected 
planar domains (cf. \cite{C, G}).

\subsection{Annulus}

 The annulus $\A_q=\{z\in \CC: q<|z|<1\}$ for $q\in (0,1)$ occupies a special 
place among multiply connected planar domains.
The first extension of the Loewner equation from
simply connected domains to  annuli goes back to Y. Komatu \cite{K1}.
Fix an annulus $\A_{Q}$ for $0<Q<1,$ and a Jordan arc $\gamma=\{\gamma(t):0\le t\le t_\gamma\}$ 
with 
$\gamma(0)\in \partial\DD$ and  $\gamma(0,t_\gamma]\subset \A_Q.$
There exists  a strictly increasing function $\alpha: [0,t_\gamma]\mapsto [Q,Q_\gamma]$
with $\alpha(t_\gamma)=Q_\gamma<1$ and
the following property: if $\alpha(t)=q,$ then there is a unique conformal map $g_q$ from $\A_Q\setminus \gamma[0,t]$ onto $\A_q$ 
such that 
$g_q(Q)=q$.  
A differential equation
for $g_q$ in the left derivative in $q$
was derived in \cite{K1} in terms of the Weierstrass 
as well as  Jacobi elliptic functions. 
But the continuity of $\alpha$ and right differentiability of $g_q$ in $q$ were not rigorously established
although an annulus variant of the Carath\'eodory convergence theorem was indicated in \cite{K1} to cover these points.

Recently \cite{FK} utilizes this variant of the Carath\'eodory theorem to show that $\alpha$ is indeed continuous and that 
 $g_q(z),\; Q\le q\le Q_\gamma,$ satisfies a genuine ODE
\begin{equation}\label{e:6.54}
\frac{\partial \log g_q(z)}{\partial\log q}
=\K_q(g_q(z),\lambda(q))-i \Im \K_q(q,\lambda(q)),\qquad  
  g_Q(z)=z,
\end{equation}
where $\K_q(z,\zeta)$, $ z\in \A_q$, $\zeta\in \partial\DD$,  is Villat's kernel defined by
$  \K_q(z,\zeta)=\K_q(z/\zeta) $.
Here 
$$ 
 \K_q(z):=\lim_{N\to\infty}\sum_{n=-N}^N \frac{1+q^{2n}z}{1-q^{2n}z},
 $$
 and $\lambda(q) 
:=g_q(\gamma(t))$ (where $t>0$ is such that $\alpha (t)=q$)
 is a continuous function of $q$ taking values on the outer boundary $\partial \DD.$
Since $\alpha$ is continuous, the curve $\gamma$ can be parametrized as $\{\gamma(q): Q\le q\le Q_\gamma\}$ so that $g_q(z)$ is a conformal map from $\A_Q\setminus \gamma[0,q]$ onto $\A_q$ 
with the 
normalization $g_q(Q)=q$. 
We may further let
$P=-\log Q$,  $P_\gamma=-\log Q_\gamma$,  $S_p(z,\zeta)=\K_{e^{-p}}(z,\zeta)$ and change the parameter $q$ into $s$ by
$q=e^{s-P}$ for $ 0\le s\le s_\gamma=P-P_\gamma.$
 Then \eqref{e:6.54}
 becomes, 
  for $z\in \A_Q\setminus \gamma[0,s]$ and $s\in [0,s_\gamma]$, 
\begin{equation}\label{e:6.55}
\frac{\partial \log  g_s(z)}{\partial s}
=S_{P-s}(g_s(z),\lambda(s))-i \Im S_{P-s}(e^{s-P},\lambda(s)),
\quad  g_0(z)=z,
\end{equation}
where  $\lambda(s):=g_s(\gamma(s))$. 
Note  that   $g_s$ is the conformal mapping
from $\A_{Q}\setminus \gamma[0,s]$ onto $\A_{Qe^s}$ with
$ g_s(Q)=Qe^s$.  
By using the stated variant of the Carath\'eodory convergence theorem,  it is also shown in \cite{FK} 
that,  the equation \eqref{e:6.54} in the right derivative sense is still valid if we take, in place of the Jordan curve $\gamma$, a family $\{F_q\}$ of growing hulls in $\A_Q$ that is right continuous with limit $\lambda(q)$ in 
  a sense similar to (5.24) of \cite{CF}.

In \cite{Z1}, D. Zhan defined    the annulus ${\rm SLE}_\kappa$   to be the growing hulls $\{K_s\}$ in $\A_Q$ 
driven by $\lambda(s)=e^{i B(\kappa s)}$  for the one-dimensional Brownian motion $B(s)$ 
based on the equation  
 \begin{equation}\label{e:7.5a}
\frac{\partial \log    g_s(z)}{\partial s}
=S_{P-s}(  g_s(z),   \lambda(s)), 
\qquad    g_0(z)=z.
\end{equation}
As was noted in the proof of \cite[Proposition 2.1]{Z1}, 
for any pair $(g_q(z), \lambda(s))$ satisfying equation \eqref{e:6.55},
its rotation $e^{i\theta(s)}(g_q(z), \lambda(s))$
satisfies \eqref{e:7.5a} where
$\theta (s)= \int_0^s \Im S_{P-r}(e^{r-P},\lambda(r))dr,$ 
and so the growing hulls based on \eqref{e:6.55} driven by $\lambda(s)$
are the same as those based on \eqref{e:7.5a} driven by $e^{i\theta(s)}\lambda(s).$

For each $\kappa >0$, 
 it was shown in \cite{Z1} 
that the distribution of the annulus ${\rm SLE}_\kappa$ 
defined by \eqref{e:7.5a} is   related to that  of  the radial ${SLE}_\kappa$ 
stopped upon hitting a compact set containing the origin.  When $\kappa=6$,  
they were further identified up to a time change.
In this connection, we point out  a gap 
in the proof of the differentiability of the function $f_t(w)$ in $t$ for each $w\in {\bf C}_0$ in \cite[page 350]{Z1}.  See Remark 2.9.
 The proof of \cite[Prop. 4.40, Th. 6.13]{L1} involves a similar gap.
         
 For a given continuous function $\lambda(q),\; Q\le q <1,$ taking value in $\partial \DD$, the ODE \eqref{e:6.54} admits a unique solution $g_q(z)
$ that can be verified to satisfy the normalization condition $g_q(Q)=q$, 
due to the fact that $\Re \K_q (q, e^{i\theta })=1$ for every $q\in (0, 1)$ and $\theta \in [0, 2\pi)$.
It may be 
worthwhile to consider an SKLE 
on the annulus based directly on the equation \eqref{e:6.54} or on its modified version
driven by 
a general diffusion process on $\partial \DD$
along the lines of \cite{CF} and this paper.

D. Zhan  further extended 
the notion of annulus ${\rm SLE}_\kappa$, $\kappa \le 4$,
 in a certain way to specify the end points of the curves and 
 investigated its 
properties such as reversibility and 
restriction property (see \cite{Z2} and references therein).  
Recently, G. Lawler \cite{L2} defined ${\rm SLE}_\kappa$ for $\kappa\le 4$ 
  in   more general multiply connected domains 
 using the Brownian loop measure and compared
   it with Zhan's one in the annulus case. 
As is noted in Remark 6.12 of \cite{CF}, 
we can hardly expect a straightforward generalization of the restriction property to the chordal ${\rm SKLE}_{\sqrt{8/3}, -b_{\rm BMD}}$ 
due to an effect of the second order BMD-domain constant $c_{\rm BMD}$.
It would be interesting to find 
connections 
of the conditional laws induced by ${\rm  SKLE}_{\sqrt{\kappa}, b}$ 
with Lawler's measures.

\subsection{Circularly slit annulus}

Parallel 
to the BMD complex Poisson kernel, the notion of the BMD Schwarz kernel $\SS(z,\zeta)$ is introduced 
in \cite{FK} for 
a general multiply connected planar domain $D$ as an analytic function in $z\in D$ whose real part is the BMD-Poisson 
kernel. In particular, it is shown in \cite{FK} that 
the Villat's kernel multiplied by $1/(2\pi)$ 
is a 
BMD Schwarz kernel for the annulus.

A domain $D$ of the form $D=\A_q\setminus \bigcup_{j=1}^{N-1} C_j$ is called a {\it circularly slit annulus} if $C_j$ are mutually disjoint concentric circular slits contained in $\A_q.$
 We denote by $\dd$ the collection of all circularly slit annuli.
We  fix $D=\A_Q\setminus \bigcup_{j=1}^{N-1} C_j \in \dd$ and
consider a Jordan arc $\gamma:[0,t_\gamma]\mapsto D$ with $\gamma(0)=\partial\DD.$
We can then find a strictly increasing function $\alpha:\; [0,t_\gamma]\mapsto [Q, Q_\gamma],\ (\alpha(t_\gamma)=Q_\gamma<1)$ such that, for $q=\alpha(t),$ there exists a unique conformal map
$g_q:\; D\setminus \gamma[0,t]\ \mapsto\ D_q=\A_q\setminus \bigcup_{j=1}^{N-1} C_j(q)\in \dd,\quad{\rm with}\quad g_q(Q)=q.
$

The first extension of the Loewner equation to a circularly slit  annulus 
goes back to Y. Komatu \cite{K2} and the resulting Komatu-Loewner equation for $g_q$ is rewritten by \cite{BF2} and then by \cite{FK} as
\begin{equation}\label{e:6.56}
\frac{\partial^- \log g_q(z)}{\partial \log q}=2\pi \wh \SS_q(g_q(z), \lambda(q)),\ q\in \alpha(0,t_\gamma]\subset (Q,Q_\gamma],\ g_Q(z)=z,
\end{equation}
where the left hand side denotes the left derivative and $\wh\SS_q(z,\zeta)=\SS_q(z,\zeta)-i\Im \SS_q(q,\zeta)$ is the normalized BMD Schwarz kernel for $D_q\in \dd$.  When $N=1$, $2\pi \wh\SS_q$ is just the normalized Villat's kernel with  the stated explicit expression and \eqref{e:6.56} is reduced to \eqref{e:6.54}.  When $N>1,$ the problem of the continuity of $\alpha$ and right differentiability of $g_q$ remains open.  Recently C. Boehm and W. Lauf \cite{BL}  establish 
a Komatu-Loewner equation for a circularly slit disk as a genuine ODE by using an extended version of the Carath\'eodory convergence theorem.
A method similar
to \cite{BL} or to \cite {CFR} might work to make \eqref{e:6.56} a genuine ODE and we may then conceive an SKLE for it
in analogue to \cite{CF}.

\subsection{Circularly slit disk} 

A domain $D$ of the form $D=\DD\setminus \bigcup_{j=1}^{N-1} C_j$ is called a {\it circularly slit disk} if $C_j$ are mutually disjoint concentric circular slits contained in  $\DD=\{z\in \CC:|z|<1\}.$  
For a circularly slit disk $D$, Bauer and  Friedrich have obtained a radial Komatu-Loewner equation \cite[(44)]{BF1} with a kernel  explicitly expressed 
in terms of the Green function and harmonic measures, which could be identical with a BMD Schwarz kernel for the  image domain $D_t$.  
The differntiability problem for $g_t(z)$ in this equation seems to have been settled by 
the aforementioned  approach of \cite{BL}.

Moreover, 
an ${\rm SKLE}_{\sqrt{\kappa},b}$ on $D$ is formulated in \cite{BF1}  
for  any 
 constant $\kappa>0$ in a quite analogous manner to \cite{CF} 
and it is  
claimed  that ${\rm SKLE}_{\sqrt{6},b}$ enjoys the locality property relative to canonical maps for a specific choice (Ansatz) of the drift coefficient $b$ of the driving process on $\partial\DD$. Our natural guess is that $b=-b_{\rm BMD}.$  However the establishment of a generalized Komatu-Loewner equation 
\cite[(63)]{BF1} for image hulls by a canonical map as a genuine ODE requires the continuity of $g_t^*$ 
(which corresponds to $\wt g_t$ in \cite{CF}). But this has been left
 unconfirmed  even in the radial SLE case with no circular slit.

{\small

\vskip 0.3truein

\noindent {\bf Zhen-Qing Chen}

\smallskip \noindent
Department of Mathematics, University of Washington, Seattle,
WA 98195, USA

\smallskip\noindent
E-mail: \texttt{zqchen@uw.edu}

\bigskip

\noindent {\bf Masatoshi Fukushima}

\smallskip \noindent
 Branch of Mathematical Science,
Osaka University, Toyonaka, Osaka 560-0043, Japan.

\smallskip\noindent
 Email: {\texttt fuku2@mx5.canvas.ne.jp}

\bigskip

\noindent {\bf Hiroyuki Suzuki}

\smallskip \noindent
Department of Physics, 
Chuo University, Tokyo 112-8551, Japan

\end{document}